\documentclass[twocolumn,journal]{IEEEtran}
\IEEEoverridecommandlockouts
\usepackage{cite}
\usepackage{amsmath,amssymb,amsfonts}
\usepackage{algorithmic}
\usepackage{graphicx}
\usepackage{textcomp}
\usepackage{amsfonts} % fraktur font
\usepackage[colorlinks, bookmarksopen, bookmarksnumbered, citecolor=blue, linkcolor=blue, urlcolor=blue]{hyperref}
\usepackage{xcolor}    % Color some contents
\usepackage{balance}   % Balance the last page.
\usepackage{booktabs}    % For table
\usepackage{amsthm}      % For citing theorem
\usepackage{floatrow}
\usepackage{subfig}
\usepackage[capitalize]{cleveref}
\def\BibTeX{{\rm B\kern-.05em{\sc i\kern-.025em b}\kern-.08em
    T\kern-.1667em\lower.7ex\hbox{E}\kern-.125emX}} 

\makeatletter
\newcommand{\linebreakand}{%
  \end{@IEEEauthorhalign}
  \hfill\mbox{}\par
  \mbox{}\hfill\begin{@IEEEauthorhalign}
}

\newtheorem{proposition}{Proposition}
\theoremstyle{definition}

\makeatother

\begin{document}
% Titles are generally capitalized except for words such as a, an, and, as at, but, by, for, in, nor, of, on, or, the, to and up, which are usually not capitalized unless they are the first or last word of the title.
% Linebreaks \\ can be used within to get better formatting as desired.
% Do not put math or special symbols in the title.
\title{3D Stochastic Geometry Model for Aerial Vehicle-Relayed Ground-Air-Satellite Connectivity
}
% A Fine-Grained Meta Distribution-Based Uplink Connectivity Analysis on Space-Air-Ground Integrated Networks
\iffalse
\author{
    \IEEEauthorblockN{Yulei Wang$^{\dagger, \ddagger}$, Yalin Liu$^{\ddagger}$, Yaru Fu$^{\ddagger}$, Yujie Qin$^{\sharp}$, and xx$^{\dagger}$} \\
    \IEEEauthorblockA{\emph{$^{\dagger}$College of Electronics and Information Engineering, South-Central Minzu University, Wuhan 430074, China}} \\
    \IEEEauthorblockA{\emph{$^{\ddagger}$School of Science and Technology, Hong Kong Metropolitan University, Hong Kong 999077, China}} \\
    \IEEEauthorblockA{\emph{$^{\sharp}$School of Information and Communication Engineering, University of Electronic Science and Technology of \\ China, Chengdu 610054, China}} \\
    \IEEEauthorblockA{E-mail: ylwang\@mail.scuec.edu.cn, \{ylliu, yfu\}@hkmu.edu.hk, yujie.qin\@kaust.edu.sa}
    \vspace{-3em}
}
\fi

% \iffalse
\author{
        \authorblockN{Yulei Wang\authorrefmark{1}\authorrefmark{2}}, 
	\authorblockN{Yalin Liu\authorrefmark{2}}, 
	\authorblockN{Yaru Fu\authorrefmark{2}}, 
	\authorblockN{Yujie Qin\authorrefmark{3}}, and        
        \authorblockN{Zhongjie Li\authorrefmark{1}}\\
	\authorblockA{
            \authorrefmark{1}College of Electronics and Information Engineering, South-Central Minzu University, Wuhan 430074, China\\
		\authorrefmark{2}School of Science and Technology, Hong Kong Metropolitan University,Hong Kong 999077, China }\\
		\authorrefmark{3}School of Information and Communication Engineering, University of Electronic Science and Technology of China, Chengdu 610054, China \\
		Email: \{ylwang, lizhongjie\}@mail.scuec.edu.cn, \{ylliu, yfu\}@hkmu.edu.hk, yujie.qin@kaust.edu.sa
	\vspace{-2.85em}
}
% \fi

\iffalse
\author{
        \authorblockN{Yulei Wang\authorrefmark{1}\authorrefmark{2}}, 
	\authorblockN{Yalin Liu\authorrefmark{2}}, 
	\authorblockN{Yaru Fu\authorrefmark{2}}, 
	\authorblockN{Yujie Qin\authorrefmark{3}}, and        
        \authorblockN{xx\authorrefmark{4}}\\
	\authorblockA{
            \authorrefmark{1}College of Electronics and Information Engineering, South-Central Minzu University, Wuhan 430074, China, Email: ylwang\@mail.scuec.edu.cn\\
		\authorrefmark{2}Hong Kong Metropolitan University, Email: \{ylliu, yfu\}\@hkmu.edu.hk }\\
		\authorrefmark{3}University of Electronic Science and Technology of China, Email: yujie.qin\@kaust.edu.sa\\
		\authorrefmark{3}xxxxxxxxx, Email: yujie.qin\@kaust.edu.sa%\\
	\vspace{-3em}
}
\fi

% \Yujie Qin \ qlzhao@must.edu.mo, \P, \S, xx$^{\natural}$, xx$^{\dagger}$, and xx$^{\S}$  lizhongjie@mail.scuec.edu.cn (Z. Li)

% If you want to put a publisher's ID mark on the page you can do it like
% this:
%\IEEEpubid{0000--0000/00\$00.00~\copyright~2015 IEEE}
% Remember, if you use this you must call \IEEEpubidadjcol in the second
% column for its text to clear the IEEEpubid mark.

% use for special paper notices
% \IEEEspecialpapernotice{(Invited Paper)}

% make the title area
\maketitle

% As a general rule, do not put math, special symbols or citations
% in the abstract or keywords.
\begin{abstract} %Space-air-ground integrated networks (GASS), which facilitate data transmission from ground users (GUs) to satellites via aerial vehicle (AV) relays, have emerged as a promising component for 6G communications. While existing research has explored the transmission performance of GASS, several key factors are often overlooked: 1) GUs and AVs are distributed on a spherical surface, rather than a flat plane, and 2) neighboring AVs are typically separated by a minimum distance to avoid collisions. To address these complexities in network modeling, this paper proposes a spherical stochastic geometry-based theoretical model to analyze the two-hop uplink transmission in GASS. Specifically, we use a homogeneous Poisson point process to model the GUs' locations and a Matérn hard-core point process for the AVs, capturing their spatial randomness. We derive the spherical distance distributions for both the GU-to-AV and AV-to-satellite transmission links, and compute their average successful transmission probabilities. Extensive Monte Carlo simulations validate the accuracy of our theoretical model, highlighting its practical value in helping practitioners configure optimal parameters for reliable GASS.
Due to their flexibility, aerial vehicles (AVs), such as unmanned aerial vehicles and airships, are widely employed as relays to assist communications between massive ground users (GUs) and satellites, forming an AV-relayed ground-air-satellite solution (GASS). In GASS, the deployment of AVs is crucial to ensure overall performance from GUs to satellites. This paper develops a stochastic geometry-based analytical model for GASS under Matérn hard-core point process (MHCPP) distributed AVs. The 3D distributions of AVs and GUs are modeled by considering their locations on spherical surfaces in the presence of high-altitude satellites. Accordingly, we derive an overall connectivity analytical model for GASS, which includes the average performance of AV-relayed two-hop transmissions. Extensive numerical results validate the accuracy of the connectivity model and provide essential insights for configuring AV deployments. 
\end{abstract}

% Note that keywords are not normally used for peerreview papers.
\vspace{-0.2cm}

\IEEEpeerreviewmaketitle

\section{Introduction}
With the proliferation of communication satellites, such as Starlink~\cite{luo2024leo}, future wireless networks are expected to provide seamless coverage for global-range applications, including remote sensing, navigation, disaster management, and other commercial uses~\cite{talgat2024maximizing}. However, due to their long-range propagation, satellites cannot support stable connections with ground users (GUs)~\cite{talgat2024maximizing}. %~\cite{talgat2024maximizing}. %qin2023energy, talgat2024stochastic
To address this issue, various emerging aerial vehicles (AVs), such as unmanned aerial vehicles and airships, can be employed as aerial relays to facilitate communication between GUs and satellites. %~\cite{liu2024space}. %qin2023energy
This \textit{AV-relayed ground-air-satellite solution (GASS)} is crucial to support data uploading of massive GUs for diverse Internet of Things tasks in the future~\cite{fu2023distributed}. In GASS, the deployment of AVs is essential to ensure overall performance from GUs to satellites. %Due to their flexibility
Particularly, multiple AVs can be widely deployed to cover %a wide range of 
massive GUs. In addition, global-range AVs can relay data to satellites to reduce the access burden from massive GU.  

%Emerging space-air-ground integrated networks (GASS) have demonstrated significant potential in providing seamless connectivity, global broadband coverage, and extensive computing services~\cite{liu2024space, qin2023energy, fu2023distributed}. GASS are expected to play a crucial role in supporting future applications such as the Internet of Things (IoT)\cite{fu2023distributed}, remote sensing\cite{tan2024outage}, and emergency communications~\cite{swaminathan2021haps} in the forthcoming 6G communications. 
%As typical two-hop networks, GASS combine the merits of terrestrial, aerial, and space networks, encompassing both ground user (GU)-to-aerial vehicle (AV) and AV-to-satellite transmission links. Specifically, GUs in remote areas collect environmental information and upload the collected data to AVs, which then relay the information to satellite acting as space-based base stations, as shown in Fig. \ref{Fig-Overview}. In this way, GASS significantly extend coverage compared to traditional terrestrial networks.

\begin{figure}[t]
    \centering
    \includegraphics[scale=0.8]{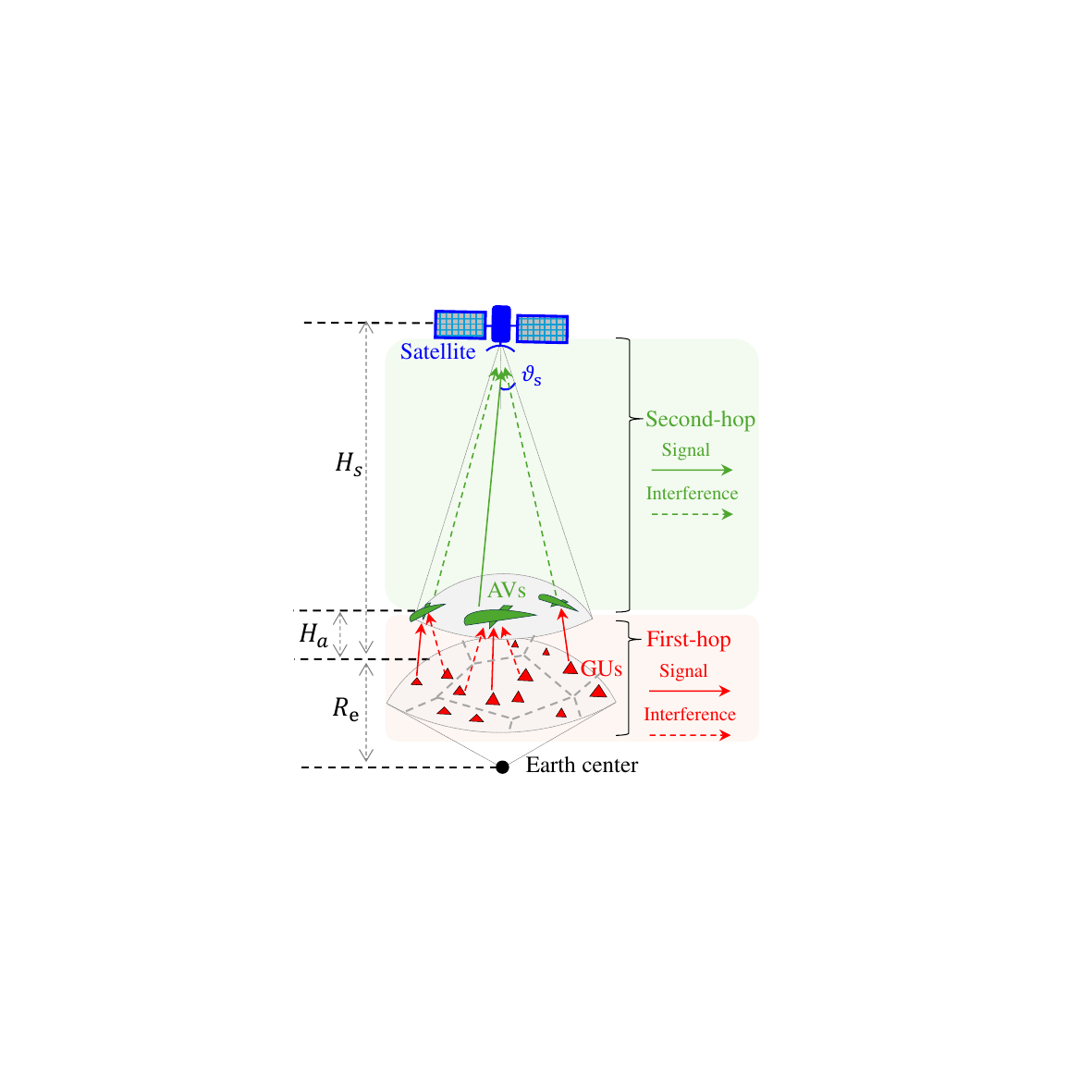}
    \caption{A typical AV-relayed two-hop transmission scenario in GASS. %ground-air-satellite solutions (GASS). 
    Herein $R_{\mathrm{e}}$ is the earth radius, $H_s$ and $H_a$ are the height of satellites and AVs over the ground, respectively.}
    \label{Fig-Overview}
    %\vspace{-5 mm} % 增加垂直空间
\end{figure}

%草稿
%stochastic geometry is a commonly used tool to model node distributions and investigate network performance metrics, e.g., successful transmission rate (also called connectivity) and outage probability~\cite{talgat2024stochastic,wang2024outage}. Based on stochastic geometry, some previous studies model AVs' distribution and investigate various network scenarios, including aerial-ground networks~\cite{liu2021ground,zhang2022multi}, non-terrestrial networks~\cite{s2a2022,ntn2024}, and space-air-ground integrated networks~\cite{uplinkGASS2024,liu2024space}.
To investigate the deployment of AVs, previous studies modeled their distributions as stochastic point processes and analyzed transmission performance among AVs, GUs, and satellites~\cite{qin2023downlink, liu2021connectivity, zhang2022multi, uplinkSAGIN2024, s2a2022, liu2024space}. Specifically, studies in~\cite{qin2023downlink, liu2021connectivity} modeled AVs' distribution using a homogeneous Poisson point process (HPPP) for %various scenarios, including 
air-ground communications and air-satellite networks, respectively. However, due to the independent nature of HPPP, random point locations in an HPPP can overlapp, it cannot cater to the practical distribution of multiple AVs, i.e., any two AVs should maintain a certain distance to avoid collisions. To address this issue, the authors in~\cite{zhang2022multi, uplinkSAGIN2024} adopted the Matérn hard-core point process (MHCPP) to model AVs' locations by introducing a minimum distance between all AVs. However, they considered the MHCPP in a plane, whereas AVs in GASS should be distributed on a 3D surface above the ground in the presence of satellites. Taking the 3D distribution into account, the works in~\cite{s2a2022,liu2024space} explore the 3D location modeling of AVs. Nevertheless,~\cite{s2a2022} mainly investigated air-satellite performance without considering the coverage of GUs. Although~\cite{liu2024space} examined the AV-relayed performance from GUs to satellites, the AVs were still modeled using an HPPP. In summary, there is a lack of performance analysis work for GASS in consideration of precise 3D MHCPP distributions for AVs.

% liu2021ground 
% \\ Ground-to-{UAV} Communication Network: Stochastic Geometry-based Performance Analysis
% zhang2022multi 
% \\ Multi-UAV Enabled Aerial-Ground Integrated Networks: A Stochastic Geometry Analysis

% s2a2022
% \\ Stochastic Analysis of Cooperative Satellite-UAV Communications
% ntn2024
% \\ System-Level Metrics for Non-Terrestrial Networks Under Stochastic Geometry Framework

% liu2024space 
% \\ Space-Air-Ground Integrated Networks: Spherical Stochastic Geometry-Based Uplink Connectivity Analysis
% uplinkGASS2024
% \\ Stochastic Geometry Based Modeling and Analysis of Uplink Cooperative Satellite-Aerial-Terrestrial Networks for Nomadic Communications With Weak Satellite Coverage

To fill in the above research gap, we aim to develop a stochastic geometry-based analytical model for GASS under 3D MHCPP-distributed AVs. However, this aim faces several challenges. First, there is no prior modeling of 3D locations for MHCPP-distributed AVs. Second, exploring the geometric relationships among AVs, served GUs, and satellites is inherently complex. Finally, considering the various transmission links among AVs, GUs, and satellites in GASS, constructing an analytical performance model is particularly challenging. To address these challenges, we explore the spherical coverage regions of satellites and AVs and model the 3D locations of AVs and GUs in the regions. Accordingly, we investigate the overall performance of GASS by finding the average transmission performance of all possible links. Overall, our key contributions are summarized as follows: 
\begin{itemize}
\item %We build a spherical geometry-based network model for GASS. 3D distributions of AVs and GUs are modeled by considering their locations on spherical surfaces because of the high-altitude satellites. Meanwhile, catering to practical node deployments, GUs in GASS are associated with their nearest AVs under an MHCPP distribution. 
We develop a %spherical geometry-based 
3D node distribution model of GASS. In the presence of high-altitude satellites, the 3D distributions of AVs and GUs are modeled by considering their locations on spherical surfaces. Catering to practical node deployments, GUs are associated with their nearest AVs; the latter is distributed as an MHCPP. 
%We adopt an HPPP and a Matérn hard-core point process (MHCPP) to model the spatial locations of GUs and AVs, respectively. These distributions effectively capture the independent scattering of GUs and the spatial repulsion characteristics of AVs in practical deployment.
\item %\textit{\textbf{Spherical Distance and Transmission Analysis}}:
We present a stochastic geometry-based analytical model for GASS. This model derives the overall connectivity for typical AV-relayed two-hop transmissions in GASS. Based on stochastic geometry, the average transmission performance of each hop is analyzed by considering all possible links from distributed transmitters of that hop. 
%We derive the spherical distance distributions for the GU-to-AV and AV-to-satellite transmission links. Based on the distributions, we express the average successful probabilities (ASPs) for these links and the overall two-hop connectivity performance of GASS. These expressions incorporate critical network parameters, such as the densities of GUs and AVs and the minimum inter-AV distance.
\item %\textit{\textbf{Validation and Practical Insights}}: 
%We validate the accuracy of the proposed theoretical model via extensive Monte Carlo simulations. The results reveal the impacts of spatial randomness on the reliability of GASS and provide valuable insights for practitioners (e.g., network engineers) to configure suitable parameters to enhance the performance of GASS.
We conduct extensive numerical analyses for the overall connectivity analytical model. The simulation results validate the accuracy of the model and provide essential insights for practical practitioners, e.g., configuring AV deployments to enhance performance in GASS. 
%The analytical results of the connectivity align with the simulation results, thereby validating the accuracy of our analytical model. Overall, our analytical model can help practitioners in estimating the practical performance of GASS in various application scenarios by adjusting system parameters. In addition, the presented analytical model can contribute many practical implementations for future studies, such as constructing objective functio
\end{itemize}

The remainder of the paper is organized as follows. Section \ref{Sec-System} introduces the system model. The 3D spherical geometry-based analytical model is presented in Section \ref{Sec-Theo}. Section \ref{Sec-Verification} gives the numerical results. Section \ref{Sec-Conclusion} concludes this paper.
%\vspace{-0.3cm}

\section{System Model}\label{Sec-System}
%\vspace{-0.3cm}

\subsection{Network Model}
We consider a typical AV-relayed two-hop transmission scenario in GASS. As shown in Fig. \ref{Fig-Overview}, a reference satellite accepts the data from multiple AVs that serve as relay nodes collecting data from numerous GUs. 
%The node distributions of AVs and GUs are modeled as follows.
%\subsubsection{Matérn hard-core point process (MHCPP) of AVs} 
%\subsubsection{MHCPP of AVs} 
Multiple AVs hover at a height $H_a$ to cover ground regions. In practice, the two neighboring AVs generally maintain a distance from each other to avoid collisions and ensure efficient coverage. In this case, we use an MHCPP of type II (called MHCPP for short) $\Phi_a$\footnote{We choose MHCPP of type II because it results in a higher density of effective points than MHCPP of type I for any density and hardcore~\cite{haenggi2011mean}.} to model the distribution of AVs and ensure a minimum distance $\check{d}$ among any two AVs. $\Phi_a$ is generated by thinning an original %homogeneous Poisson point process (HPPP) 
HPPP $\Phi_a^0$. Let $\lambda_a^0$ be the node density in $\Phi_a^0$ and $p_a\left(\check{d} \right)$ be the probability of retaining an AV in $\Phi_a^0$. $p_a\left(\check{d} \right)$ is given by~\cite{wang2024analytical}
\begin{equation}
p_a\left(\check{d} \right) = \left\{\begin{matrix}\frac{1 - \exp\left(-\lambda_a^0\pi{\check{d}}^2\right)}{\lambda_a^0\pi{\check{d}}^2}&\check{d}>0\\1&\check{d}=0\\\end{matrix}.\right.	
\label{Eq-pdfNakagami}
\end{equation}
After the thinning, $\Phi_a$ has the node density $\lambda_a = \lambda_a^0 p_a\left(\check{d} \right)$. %The MHCPP $\Phi_a$ has identical distribution characteristics with an HPPP $\Tilde{\Phi}_a$ with the same density $\lambda_a$, as validated in~\cite{haenggi2011mean}. %, i.e., $\Tilde{\Phi}_a \approx \Phi_a$, and the accuracy of such approximation is validated in \cite{haenggi2011mean}.
The GUs are uniformly scattered on the Earth's surface (i.e., ground). Hence, we model the location distribution of GUs as an HPPP $\Phi_u$ with density $\lambda_u$. Each GU associates with its geographically nearest AV for transmission, where a Voronoi tessellation of $\Phi_u$ according to multiple AVs is formed~\cite{wang2024outage}. Meanwhile, each GU is considered to have data to transmit to its associated AV with the probability $p_u^t$. Using the thinning process of HPPP, the transmitting GUs form an independent HPPP $\Phi_u^t$ ($\Phi_u^t \subset \Phi_u$) with density $\lambda_u^t = p_u^t \lambda_u$. In the presence of satellites, the 3D locations of AVs and GUs are modeled based on a spherical surface, as presented in Section~\ref{subsec: 3Ddistribution}.

%\vspace{-0.4cm}
\subsection{Transmission Model}
\label{SubSec-Trans}
The two-hop transmission in GASS includes: i) the first-hop link from a reference GU to its associated AV, and ii) the second-hop link from the AV to a reference satellite. %According to Slivnyak’s theorem \cite{xx}, it focuses on a typical GU and AV to evaluate the overall performance of the network. 
For convenience, the subscript $i, i \in \{1, 2\}$ is adopted to indicate the first link and the second link, respectively.
%\subsection{Signal Propagation Model}\label{SubSec-Propa}
%The signals in two-hop uplink transmissions %from GUs to the satellite with the relay of the AVs 
%undergo attenuation of various effects, such as antenna gains, small-scale and large-scale fading. %Below, we illustrate these effects, respectively.
%\subsubsection{Antenna Gains}
%Antenna design and direction influence effective antenna gain. 
Let $G_i$ denote the receiver antenna gain of the $i$-th hop link. Its value is given by $G_i = \iota_i \left(\pi D_i f_i/c\right)^2$, where $\iota_i$ is the antenna illumination coefficient, $f_i$ is the carrier frequency used, $D_i$ is the diameter of the normalized reflector antenna, and $c$ is the light speed~\cite{liu2024space}. 
%\subsubsection{Channel Fading}
For both links $i \in \{1, 2\}$, line-of-sight (LoS) propagation dominates the small-scale fading process. We use Nakagami-m fading to capture the channel fading characteristics. Let $h_i$ be the channel power gain of the $i$-th hop link, %. $h_i$ is the square of the amplitude fading coefficient and %follows the gamma distribution. Hence, 
its probability density function (PDF) is given by $\frac{m_i^{m_i}}{\Gamma(m_i)\Omega_i^{m_i}} h_i^{m_i-1} \exp\left(-\frac{m_i}{\Omega_i} h_i\right)$~\cite{song2022cooperative},
% \begin{equation}
% \begin{aligned}
% f(h_i; m_i, \Omega_i) = \frac{m_i^{m_i}}{\Gamma(m_i)\Omega_i^{m_i}} h_i^{m_i-1} \exp\left(-\frac{m_i}{\Omega_i} h_i\right),
% \end{aligned}
% \label{Eq-pdfNakagami}
% \end{equation}
where $m_i$ is the shape parameter, $\Omega_i$ is the spread parameter, and $\Gamma(\cdot)$ is the standard gamma function \cite{sklar2021digital}.
%Path loss of all transmission links can be approximated to the free-space model. 
Let $r_i(\texttt{t})$ be the distance between the transmitter $\texttt{t}$ and its receiver in the $i$-th hop link, the corresponding path loss $L_i(\texttt{t})$ %
%Let $L_i$ be the path loss of a $i$-th hop link, which 
is given by~\cite{jung2022performance}
\begin{equation}
\begin{aligned}
L_i(\texttt{t}) = l_i \left(\frac{c}{4\pi f_i}\right)^2 (r_i(\texttt{t}))^{-2},%{-\alpha},
\end{aligned}
\label{Eq-pathloss}
\end{equation}
where $l_i$ is the additional loss caused by atmospheric effects and rain/fog attenuation. %, and $\alpha$ is the path loss exponent. % (i.e., $\alpha = 2$ in free space).
%\vspace{-0.4cm}
\subsection{Interference Model}
Generally, frequency division multiple access technology is adopted to enable simultaneous transmissions from multiple transmitters (i.e., GUs and AVs) in each hop. However, the available carriers are limited compared to the number of transmitters under the wide coverage of GASS. Each link suffers interference when multiple transmitters choose the same carrier to transmit. To evaluate the link quality, the signal-to-interference-plus-noise ratio (SINR) at the receiver is measured.
%Given a carrier, l
Let $\texttt{SINR}_i$ be the SINR of the $i$-th hop link from a transmitter $\texttt{t}_0$ to its receiver, %we have %$\texttt{SINR}_i(\texttt{t}_0)$ 
which is given by
\begin{equation}
\begin{aligned}
\texttt{SINR}_i = \frac{P_i G_i h_i L_i(\texttt{t}_0)}{I_i + W_i},
\end{aligned}
\label{Eq-SINR}
\end{equation}
where $P_i$ is the transmission power of transmitters in the $i$-th hop link and %. %$h_i$ is the channel power gain and the distance of the typical transceiver pair. In (\ref{Eq-SINR}), 
$I_i = \sum_{\texttt{t} \in \Phi_i^\prime}{P_i G_i h_i L_i(\texttt{t})}$ is the aggregated power of interfering signals at the receiver, where $\Phi_i^\prime=\Phi_i\setminus\{\texttt{t}_0\}$ is the set of interfering transmitters with 
$\Phi_1 = \Phi_u^t, \Phi_2 = \Tilde{\Phi}_a$.
% For $i=\{1,2\}$, $\Phi_i=\{\Phi_u^t,\Tilde{\Phi}_a\}$. %GUs and AVs, and $\Phi_1^\prime = \Phi_u^t / \{\mathrm{typical \ transmitting \ GU} \}, \Phi_2^\prime = \Tilde{\Phi}_a / \{ \mathrm{typical \ transmitting \ AV} \}$.
%Besides, 
$W_i$ is an additive white Gaussian noise. % with zero mean and variance $\sigma_i^2$.
%\vspace{-0.3cm}
\section{3D Stochastic Geometry %-Based Connectivity 
Model}
\label{Sec-Theo}
This section analyzes the two-hop transmission performance in GASS, which is greatly affected by the interference in two-hop links. To characterize the interference, we first model the %stochastic 
3D distance distributions among GUs, AVs, and satellites for each hop. Based on stochastic geometry of 3D node distributions, we investigate the overall connectivity %of two-hop transmissions in
of GASS.
%\vspace{-0.6cm}

\subsection{3D Distribution Model}%待定
\label{subsec: 3Ddistribution}
Let $\mathcal{A}_i(i \in \{1, 2\})$ be the coverage region of a receiver for its transmitters in the $i$-th hop link. The ground is approximated as a spherical surface with the Earth center $\mathbf{O}$. Given all transmitters at the same altitude, $\mathcal{A}_i$ can be approximated as a spherical cap. For the $i$-th hop link, all transmitters under the same receiver have the potential to transmit their data. Due to the distribution of transmitters, the distances $r_i(\texttt{t})$ from different transmitters $\texttt{t}$ to the same receiver are different. Below, we evaluate the distribution of $r_i(\texttt{t})$ for two-hop links.

\subsubsection{$\mathcal{A}_1$ and $r_1(\texttt{t})$}
For the first hop link, given a reference AV, $\mathcal{A}_1$ is the AV's coverage for GUs and $r_1(\texttt{t})$ is the distance between the AV and any one of GUs $\texttt{t}$ under its coverage. All AVs are generally configured with adaptive antennas to serve GUs dynamically from all directions~\cite{alshbatat2010performance}. Recall that the GUs and their associated AVs form a Voronoi tessellation; the associated coverage of each AV forms a Voronoi cell. The mean area of these Voronoi cells is given by $1/\lambda_a$~\cite{wang2024outage}, which can be approximated to a spherical cap $\mathcal{A}_1$. %~\cite{chen2018qos,liu2020analysis}
Referring to Archimedes’ hat theorem, the cap $\mathcal{A}_1$'s height $H_1$ is given by $H_1 = \left|\mathcal{A}_1\right|/(2\pi R_{\mathrm{e}}) = 1/(2\pi\lambda_a R_{\mathrm{e}})$. According to the geometry relationship in Fig. \ref{Fig-Overview}, the minimum and maximum distances of $r_1(\texttt{t})$ are given by $\check{R}_1,\hat{R}_1$ as follows.
\begin{align}
\check{R}_1 = H_a,\hat{R}_1 =\sqrt{H_a^2 + 2 H_1 \left(R_{\mathrm{e}} + H_a\right)}.
\label{eq: r1_min_max}
\end{align}
%$\check{R}_1 = H_a$ and $\hat{R}_1 =\sqrt{H_a^2 + 2 H_1 \left(R_{\mathrm{e}} + H_a\right)}$, respectively. 
Since the GUs follow an HPPP, the PDF $f(r_1)$ of $r_1(\texttt{t})$ on $\mathcal{A}_1$ can be derived by
\begin{equation}
\begin{aligned}
f(r_1)
& %= \frac{\mathbb{P}\left(x \leq r_1\right) }{\mathrm{d} r_1} 
= \frac{\mathrm{d}\frac{\left|\mathcal{A}_1(r_1)\right|}{\left|\mathcal{A}_1(\hat{R}_1)\right|}}{\mathrm{d} r_1} 
= \frac{\frac{\pi \lambda_a R_{\mathrm{e}} \left(r_1^2 - \check{R}_1^2\right)}{\left(R_{\mathrm{e}} + \check{R}_1\right)}}{\mathrm{d} r_1} 
= \frac{2\pi \lambda_a R_{\mathrm{e}} r_1}{R_{\mathrm{e}} + \check{R}_1}.
\label{Eq-PDF-R1}
\end{aligned}
\end{equation}%
% \begin{equation}
% \begin{aligned}
% & F_{r_1}\left(r\right) 
% = \resizebox{0.75\hsize}{!}{$ \mathbb{P}\left(r_1 \leq r\right) 
% = \frac{\left|\mathcal{A}_1(r)\right|}{\left|\mathcal{A}_1(\hat{R}_1)\right|} 
% = \frac{\pi \lambda_a R_{\mathrm{e}} \left(r^2 - \check{R}_1^2\right)}{\left(R_{\mathrm{e}} + \check{R}_1\right)}, $}
% \label{Eq-CDF-R1}
% \end{aligned}
% \end{equation}
%where $H_1$ is the cap $\mathcal{A}_1$'s height and $r\in[\check{R}_1,{\hat{R}}_1]$ with $\check{R}_1,{\hat{R}}_1$ being the minimum and maximum distance of $r_1$. According to Archimedes’ hat theorem, $H_1 = \left|\mathcal{A}_1\right|/(2\pi R_{\mathrm{e}}) = 1/(2\pi\lambda_a R_{\mathrm{e}})$. Following the geometrical relationship in Fig. \ref{Fig-Overview}, $\check{R}_1 = H_a$ and $\hat{R}_1 =\sqrt{H_a^2 + 2 H_1 \left(R_{\mathrm{e}} + H_a\right)}$. Referring to~\eqref{Eq-CDF-R1}, the probability density function (PDF) of $r_1$ is derived by
% \begin{equation}
% \begin{aligned}
% & f_{r_1}\left(r\right) 
% = \frac{\mathrm{d} F_{r_1}\left(r\right)}{\mathrm{d} r} 
% = \frac{2\pi \lambda_a R_{\mathrm{e}} r}{R_{\mathrm{e}} + \check{R}_1}.
% \label{Eq-PDF-R1}
% \end{aligned}
% \end{equation}

\subsubsection{$\mathcal{A}_2$ and $r_2(\texttt{t})$}
For the second hop link, given a reference satellite, $\mathcal{A}_2$ is the satellite's coverage for AVs and $r_2(\texttt{t})$ is the distance between the satellite and any one of AVs $\texttt{t}$ under its coverage. The satellite is considered to cover the earth surface vertically with its antenna. Let $\vartheta_s$ be the $3$dB beamwidth of the satellite, which is given by $\vartheta_s = c \kappa_s/(f_2 D_2)(\mathrm{degrees})$. 
%According to the geometry relationship in Fig. \ref{Fig-Overview}
Similar to~\eqref{eq: r1_min_max}, the minimum and maximum distances of $r_2(\texttt{t})$ are given by $\check{R}_2,\hat{R}_2$ as follows.
\begin{align}
&\check{R}_2 = H_s - H_a, \hat{R}_2  =(R_{\mathrm{e}} + H_s)\cos{\vartheta_s}\notag\\
&-\sqrt{(R_{\mathrm{e}} + H_a)^2 -(R_{\mathrm{e}} + H_s)^2\sin^2{\vartheta_s}}.
\end{align}
% As shown in Fig. \ref{Fig-Overview}, the minimum and maximum distance of $R_2$ can be expressed as $\check{R}_2 = H_s - H_a$ and $\hat{R}_2 =(R_{\mathrm{e}} + H_s)\cos{\vartheta_s}-\sqrt{(R_{\mathrm{e}} + H_a)^2 -（R_{\mathrm{e}} + H_s)^2\sin^2{\vartheta_s}}$, respectively.
% % \begin{equation}
% % \begin{aligned}
% %     & \check{R}_2 = H_s - H_a, \\
% %     & \hat{R}_2 = \resizebox{0.9\hsize}{!}{$ \left(R_{\mathrm{e}} + H_s\right)\cos{\vartheta_s}-\sqrt{\left(R_{\mathrm{e}} + H_a\right)^2 -\left(R_{\mathrm{e}} + H_s\right)^2\sin^2{\vartheta_s}}. $} % \resizebox{0.85\hsize}{!}{$ 在这里打公式内容 $}
% % \end{aligned}
% % \end{equation}
Referring to Archimedes’ hat theorem, the cap $\mathcal{A}_2$'s height and the area of $\mathcal{A}_2$ are given by $H_2 = (\hat{R}_2^2 - \check{R}_2^2)/(2\left(R_{\mathrm{e}} + H_s\right))$ and $\left|\mathcal{A}_2\right| = 2 \pi \left(R_{\mathrm{e}} + H_a\right) H_2$, respectively.
% \begin{equation}
% \begin{aligned}
%     \left|\mathcal{A}_2\right| 
%     & = 2 \pi \left(R_{\mathrm{e}} + H_a\right) H_2.
%     % & = \frac{\pi \left(R_{\mathrm{e}} + H_a\right) (\left(r\right)^2 - \left(H_s - H_a\right)^2)}{\left(R_{\mathrm{e}}^i + \check{R}_s^i\right)}, r_s^i\in\left[\check{R}_s^i,{\hat{R}}_s^i\right]
% \end{aligned}
% \end{equation}
Recall that the AVs follow an MHCPP (type II) $\Phi_a$. For interference analysis of multiple AVs, $\Phi_a$ can be approximated by an HPPP $\Tilde{\Phi}_a$ with the same density $\lambda_a$, as validated in~\cite{haenggi2011mean}. Therefore, we consider the approximate distribution $\Tilde{\Phi}_a$ for AVs. 
%Recall that the AVs are uniformly scattered on the AV-flying plane following an HPPP $\Tilde{\Phi}_a$. 
Considering GUs $\Tilde{\Phi}_a$, the PDF of 
%Let $f_{R_2}(r)$ be the probability density function (PDF) of $R_2$, where $r\in[\check{R}_2,{\hat{R}}_2]$. Similar to~\eqref{Eq-PDF-R1}, $f(r_2)$  of 
$r_2(\texttt{t})$ on $\mathcal{A}_2$ can be derived by
% \begin{equation}
% \begin{aligned}
% \resizebox{0.75\hsize}{!}{$
% F_{R_2}\left(r\right) = \frac{r^2 - \check{R}_2^2}{\hat{R}_2^2 - \check{R}_2^2}, \
% f_{R_2}\left(r\right) = \frac{2 r}{\hat{R}_2^2 - \check{R}_2^2}, $}
% \label{Eq-FRd}
% \end{aligned}
% \end{equation}
\begin{equation}
\begin{aligned}
f(r_2)
& %= \frac{\mathbb{P}\left(R_2 \leq r\right) }{\mathrm{d} r} 
= \frac{\mathrm{d}\frac{r_2^2 - \check{R}_2^2}{\hat{R}_2^2 - \check{R}_2^2}}{\mathrm{d} r_2} 
= \frac{2 r_2}{\hat{R}_2^2 - \check{R}_2^2}.
\label{Eq-PDF-R2}
\end{aligned}
\end{equation}%

\subsection{Overall Connectivity Model}%待定$\overline{\mathcal{P}}_{\mathrm{overall}}$
We investigate the overall connectivity of two-hop transmissions in GASS by including the average successful probability (ASP) of two-hop transmissions. Based on the stochastic geometry, the ASP of each hop is evaluated as the average transmission performance of all possible links in that hop. Let $\overline{\mathcal{P}}_{\mathrm{overall}},\overline{\mathcal{P}}_i(i\in\{1,2\})$ be the overall connectivity and ASPs of two-hop links. %, respectively. 
%We construct the formula of 
The performance $\overline{\mathcal{P}}_{\mathrm{overall}}$ is evaluated by% follows
\begin{align}
    &\overline{\mathcal{P}}_{\mathrm{overall}} \triangleq \prod_{i\in\{1,2\}}\overline{\mathcal{P}}_i,%\label{Eq-Overall}\\
    %&
    \overline{\mathcal{P}}_i %\triangleq  \mathbb{E}_{\Phi_i}\left[\mathcal{P}_i\right] 
    = \int_{\check{R}_i} ^{\hat{R}_i}{ \mathcal{P}_i f\left(r_i\right)\mathrm{d} r_i},
\label{Eq-Mean-CSP}
\end{align}%
where $\mathcal{P}_i$ is the conditional successful probability of a particular $i$-th hop link. %, $f(r_i)$ is the PDF of $R_i$, and $r_i \in [\check{R}_i, \hat{R}_i]$. 
$\mathcal{P}_i$ represents the probability of the signal from the transmitter to the receiver being successfully received and decoded, i.e., the SINR measured at the receiver exceeds a certain threshold. % given a realization of transmitters. 
Let $\theta_i$ be the SINR threshold of transmitters in the $i$-th hop link. 
%Let $\mathcal{P}_i\left(\theta_i \right)$ denote the conditional success probability (CSP) of $i$-th hop transmission link. Given a SINR threshold $\theta_i$, conditioned on $\Phi_i$ (i.e., $\Phi_1 = \Phi_u, \Phi_2 = \Tilde{\Phi}_a$), 
$\mathcal{P}_i\left(\theta_i\right)$ is given by
\begin{equation}
\begin{aligned}
     & \mathcal{P}_i\triangleq \mathbb{P}\left(\texttt{SINR}_i >\theta_i \middle|\Phi_i\right).
\label{Eq-CSP}
\end{aligned}
\end{equation}
%where $\Phi_1$, $\Phi_2$ is the SINR threshold for the signal reception at the AVs and satellite, respectively.
% The CSP captures the successful probability of individual transmission link. In order to characterize the average performance of all links in GASS, we define the average successful probability $\bar{\mathcal{P}}_s\left(\theta\right)$, i.e., the mean of CSP, as:
% \begin{equation}
% \begin{aligned}
%     \overline{\mathcal{P}}_i\left(\theta\right) & \triangleq  \mathbb{E}_{\Phi_i}\left[\mathcal{P}_i\left(\theta_i \right)\right] = \int_{\check{R}_i} ^{\hat{R}_i}{ \mathcal{P}_i\left(\theta_i \right) f(r_i)\mathrm{d} r_i},
% \label{Eq-Mean-CSP}
% \end{aligned}
% \end{equation}
% where $f(r_i)$ is the PDF of $R_i$, $r_i \in [\check{R}_i, \hat{R}_i]$. 
%Below, we derive the expression of CSP $\mathcal{P}_i\left(\theta_i\right)$.
%\subsubsection{$\mathcal{P}_i$}
Substituting %Eqs.~\eqref{Eq-SINR},~\eqref{Eq-PDF-R1}, and~\eqref{Eq-PDF-R2} 
Eq.~\eqref{Eq-SINR} into~\eqref{Eq-CSP}, we have~\cref{pro1}.
\begin{proposition}
\label{pro1}
$\forall i\in \{1,2\}$, $\mathcal{P}_i$ is given by
%\vspace{-0.4cm}

{\small
\begin{align*}
& \mathcal{P}_{i} = \exp\left(-\dot{S}_{i} - \dot{R}_{i} \varepsilon_i \right) \\
& \times \sum_{k=0}^{m_i-1} \sum_{l,q,k} \left[\frac{1}{i!j!...q!} \left(\dot{S}_i + \dot{R}_i \varepsilon_i'\right )^i\left (\dot{R}_i\varepsilon_i'' \right )^j...\left (\dot{R}_i \varepsilon_i ^{(l)} \right)^q \right],\\
& \text{where } 
\resizebox{0.9\hsize}{!}{$\dot{S}_i =  \frac{16 m_i \theta_i \sigma_i^2 (r_i(\texttt{t}_0))^2}{\Omega_i P_i \iota_i l_i \mathrm{D}_i^2 },
\dot{R}_i = \left\{\frac{\pi \lambda_i R_{\mathrm{e}}}{R_{\mathrm{e}} + H_a}, \frac{ \pi \lambda_i (R_{\mathrm{e}} + H_a)}{R_{\mathrm{e}} + H_s}\right\}$},\\
&\varepsilon_i = \int_{\check{R}_i^2}^{\hat{R}_i^2}{ \left(1- \left(1 + \frac{\theta_i (r_i(\texttt{t}_0))^2 }{\gamma} \right)^{-m_i} \right) \mathrm{d}\gamma},\gamma=(r_i(\texttt{t}))^2,\\%(r_i(\texttt{t}))^2 \in \left[\check{R}_i^2, \hat{R}_i^2\right],\\
&\resizebox{0.95\hsize}{!}{$\dot{\varepsilon}_i^{(l)} = C_l^{m_i+l-1} \int_{\check{R}_i^2 } ^{\hat{R}_i^2 } \left( \frac{\theta_i (r_i(\texttt{t}_0))^2 }{\gamma} \right)^{l} \left(1 + \frac{\theta_i (r_i(\texttt{t}_0))^2 }{\gamma}\right)^{-m_i-l}\mathrm{d}\gamma$},\\
& %\iota_i = \{\iota_a, \iota_s \}, %l_i = \{l_1, l_2 \}, 
%D_i = \{D_a, D_s \},%\lambda_i = \{\lambda_u^t, \lambda_a \},\\
\lambda_i = \{p_u^t \lambda_u, \lambda_a^0 p_a\left(\check{d} \right) \},\check{R}_i = \{H_a,H_s-H_a\},
%\label{Eq-CSP-Detailed}
\end{align*}
}%
and {\small $\hat{R}_i =\{ \sqrt{H_a^2 + 2 H_1 \left(R_{\mathrm{e}} + H_a\right)},(R_{\mathrm{e}} + H_s)\cos{\vartheta_s}-\sqrt{(R_{\mathrm{e}} + H_a)^2 -(R_{\mathrm{e}} + H_s)^2\sin^2{\vartheta_s}}\}$}.
\end{proposition}
\textit{Proof}: The proof is given in Appendix~A.
\hfill$\blacksquare$
\\

\textbf{Remarks:} Substituting~\cref{pro1} into (\ref{Eq-Mean-CSP}), we can calculate $\overline{\mathcal{P}}_{\mathrm{overall}}$. Referring to~\cref{pro1} and (\ref{Eq-Mean-CSP}), it can be seen that our analytical model can be used to investigate the effect of comprehensive system parameters, including two-hop topology settings $\{H_a,H_s,R_{\mathrm{e}}\}$, transceiver configurations $\{P_i,f_i,p_u^t,\iota_i, D_i,\kappa_s,\theta_i, r_i(\texttt{t}_0)\}$, distribution parameters $\{\lambda_u, \lambda_a^0, \check{d}\}$, and channel factors $\{m_i,\Omega_i,l_i\}$. To analyze the impact of these parameters on our model, we conduct extensive numerical results in Section~\ref{Sec-Verification}. Due to the length limit, several significant parameters $\{\check{d},\lambda_a^0,\lambda_u,\theta_1,\theta_2\}$ are analyzed.
%\vspace{-0.4cm}

\begin{figure}[tbp]
    \centering
    %\subfloat[$\overline{\mathcal{P}}_{success}$ VS. $\lambda_a^0$]
    \subfloat[$\lambda_u$ = 50GUs/Km$^2$, $\theta_1$ = 0dB, $\theta_2$ = -5dB.]
    {\includegraphics[width=7.97cm]{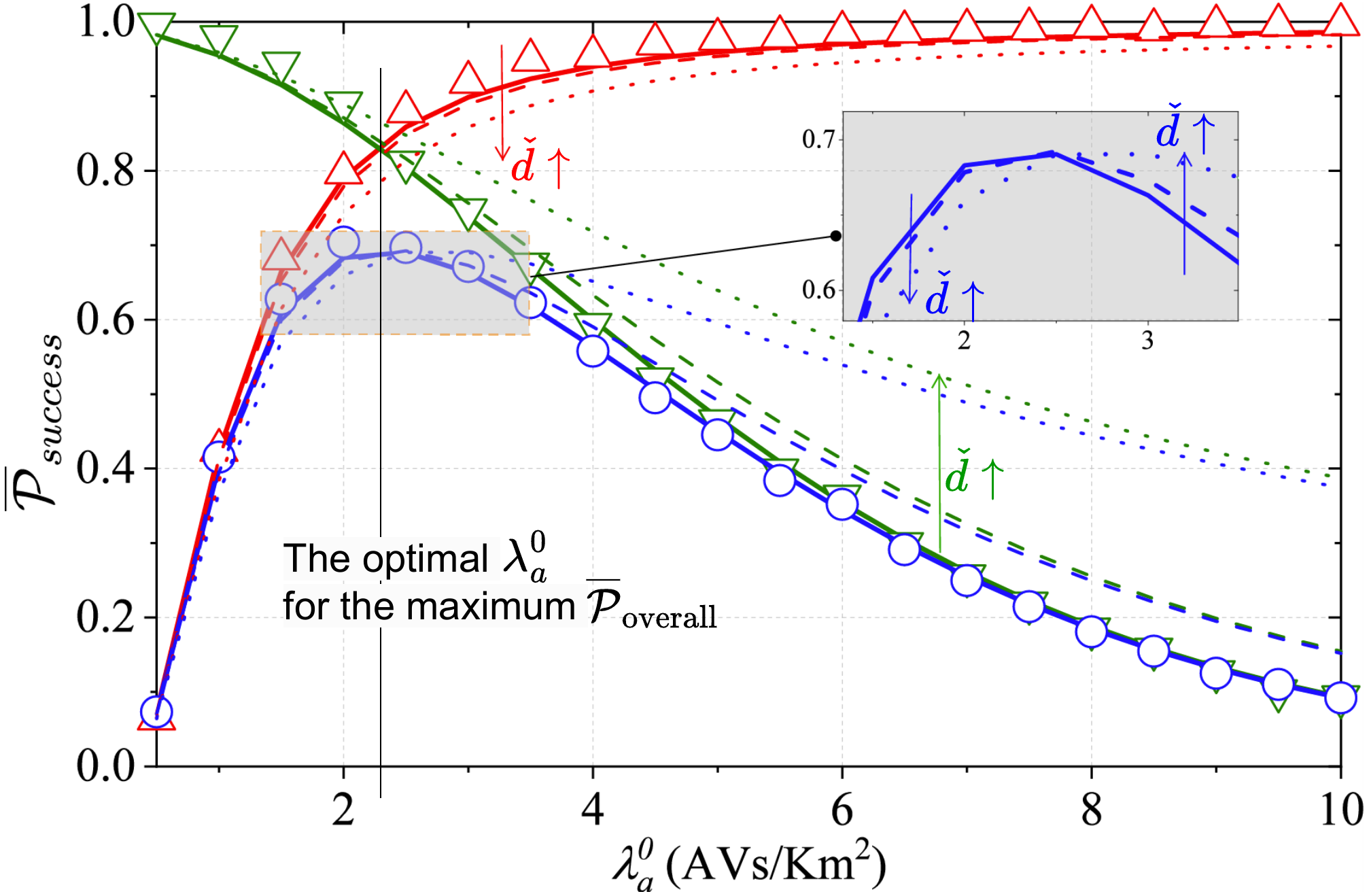}\label{Fig-Ps-lambda-a0}}	  \hfil 
    %\\[-0.2ex]
    %\subfloat[$\overline{\mathcal{P}}_{success}$ VS. $\lambda_u$]
    \subfloat[$\lambda_a^0$ = 5AVs/Km$^2$, $\theta_1$ = 0dB, $\theta_2$ = -5dB.]
    {\includegraphics[width=8cm]{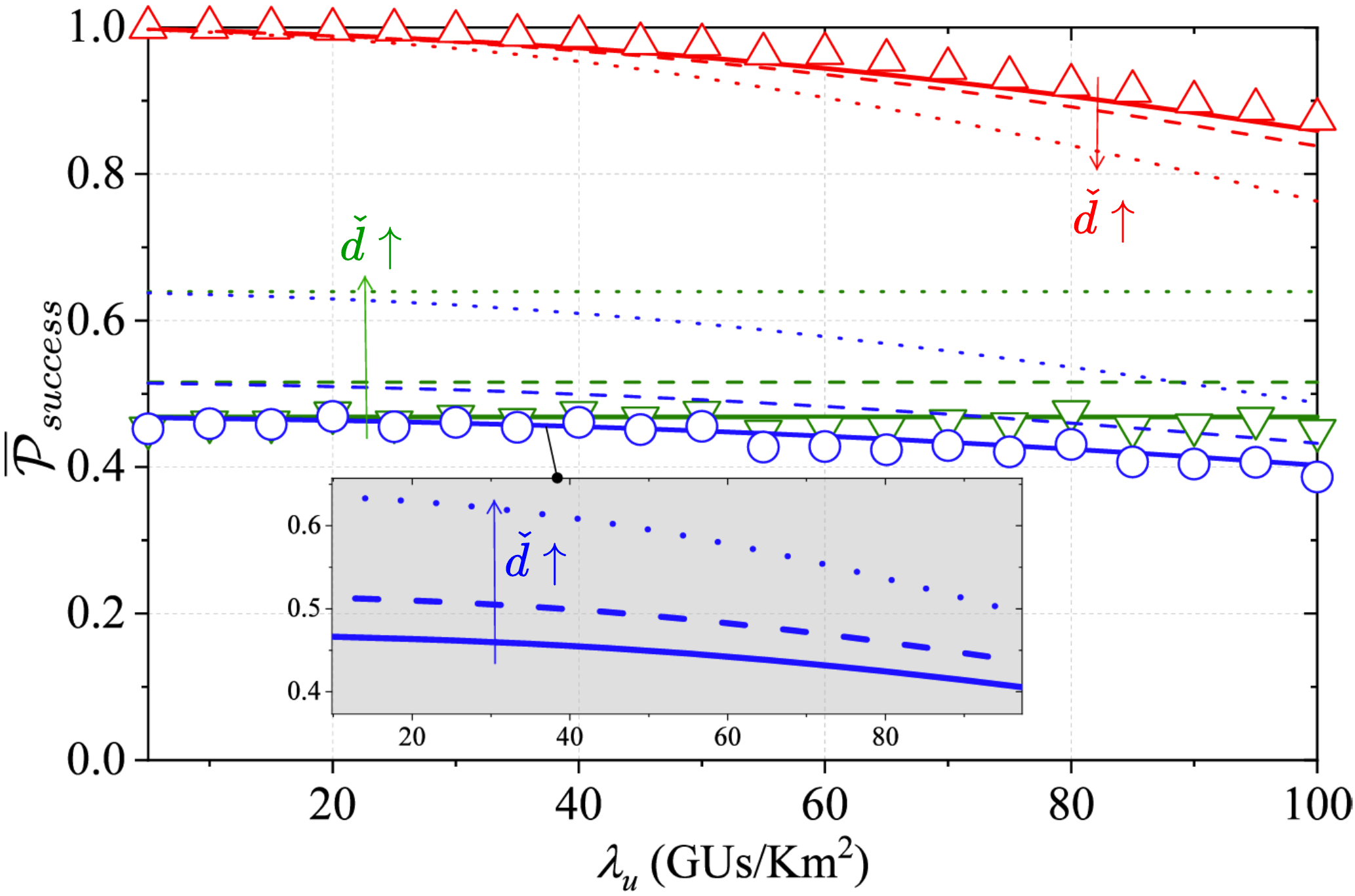}\label{Fig-Ps-lambda-u}}	     \hfil 
    %\\[-0.2ex]
    %\subfloat[$\overline{\mathcal{P}}_{success}$ VS. $\theta_i$]
    \subfloat[$\lambda_a^0$ = 5AVs/Km$^2$, $\lambda_u$ = 50GUs/Km$^2.$]%, $\theta_1$ = $\theta_2$.]
    {\includegraphics[width=7.97cm]{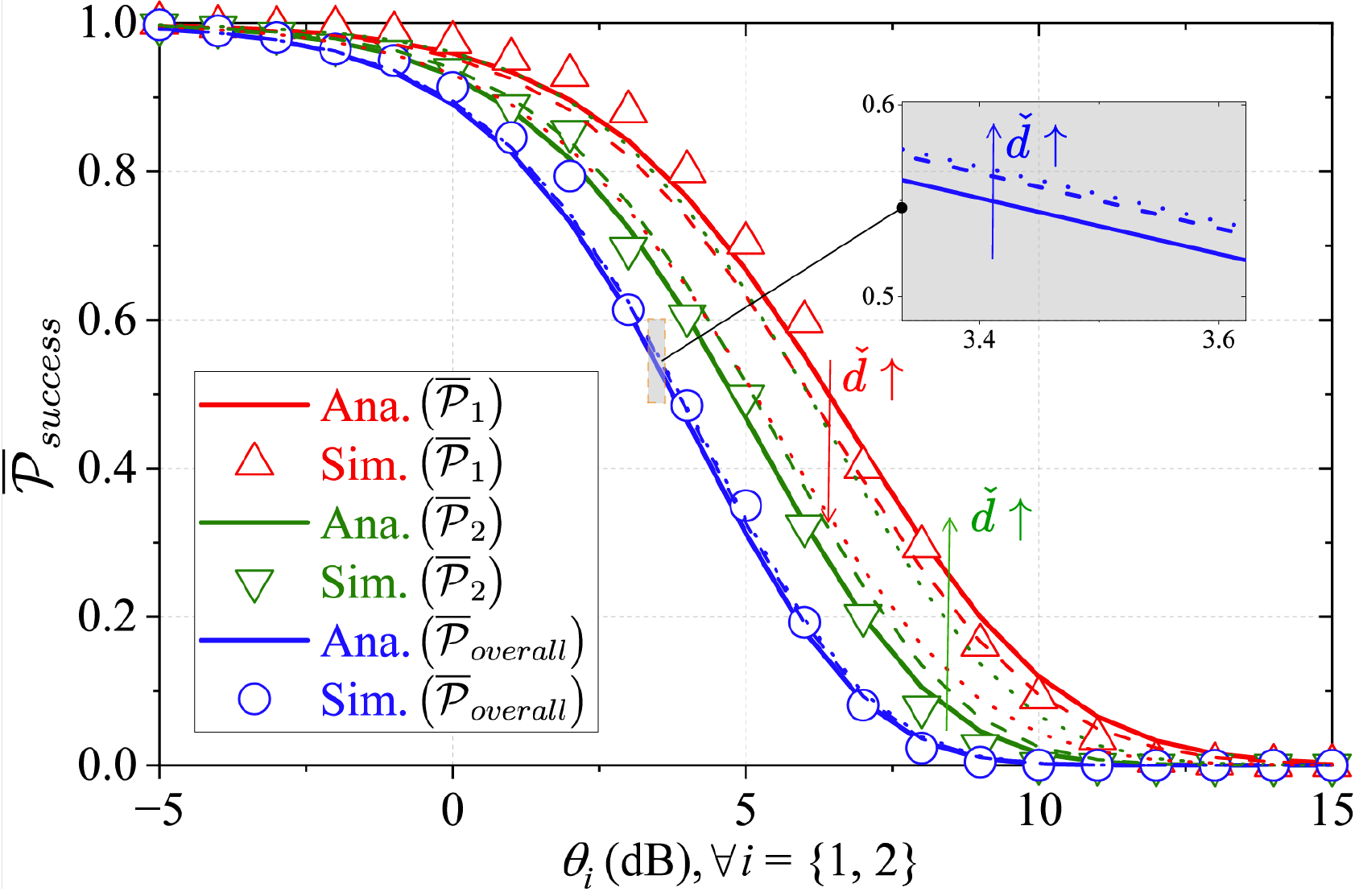}\label{Fig-Ps-theta}}	     \hfil 
    \caption{Three metrics of $\overline{\mathcal{P}}_{success}:\overline{\mathcal{P}}_1,\overline{\mathcal{P}}_2,\overline{\mathcal{P}}_{\mathrm{overall}}$ versus different parameters $\{\check{d},\lambda_a^0,\lambda_u,\theta_1,\theta_2\}$, where $H_s$ = 600Km, $H_a$ = 1Km, $\{m_1,m_2\}=\{3,3\}$, and $\{\Omega_1,\Omega_2\}=\{1,1\}$. %In each subfigure, t
    The ``solid", ``dash", and ``dot" lines %in each color group of ``red", ``green", and ``blue" 
    show the theoretical results under $\check{d}=\{0,100,200\}$m respectively. The simulation results under $\check{d} = \{100,200\}$m are not shown to avoid overlapping.} 
    \label{Fig-Ps}
    %\vspace{-0.3cm}
\end{figure}

\section{Numerical Results}\label{Sec-Verification}
This section presents the numerical results of three metrics $\overline{\mathcal{P}}_1,\overline{\mathcal{P}}_2,\overline{\mathcal{P}}_{\mathrm{overall}}$ as given in Secion~\ref{Sec-Theo}. %validates the accuracy of the proposed theoretical model through extensive Monte Carlo simulations. 
%The simulation environment is built upon the system model detailed in Section \ref{Sec-System}. 
The parameter settings are set aligned with our previous work \cite{liu2024space}. We conduct Monte Carlo simulations %The simulation results 
and take the average results over a total of ${10}^4$ iterations. In each simulation, the random distributions of GUs and AVs as well as the random channel fading of both links are generated. In our resulting Fig.~\ref{Fig-Ps}, the labels ``Ana.'' and ``Sim.'' represent the theoretical and simulation results, respectively. Besides, one legend remains for all subfigures to ensure clarity of the plotted results and avoid overlapping.

%\vspace{-0 mm}
Fig.~\ref{Fig-Ps}\subref{Fig-Ps-lambda-a0} shows the impact of $\lambda_a^0,\check{d}$ on $\overline{\mathcal{P}}_1, \overline{\mathcal{P}}_2, 
\overline{\mathcal{P}}_{\mathrm{overall}}$. We can see that the increasing $\lambda_a^0$ leads to the increasing $\overline{\mathcal{P}}_1$ and decreasing $\overline{\mathcal{P}}_2$; the larger $\check{d}$ leads to the smaller $\overline{\mathcal{P}}_1$ and the larger $\overline{\mathcal{P}}_2$. This is because a larger $\lambda_a^0$ means more relays (i.e, AVs), thus reducing the coverage burden of each AV and enhancing $\overline{\mathcal{P}}_1$. More AVs also introduce greater interference in the second hop, then lowering $\overline{\mathcal{P}}_2$. Meanwhile, the larger $\check{d}$ can reduce the AVs density, thus increasing the coverage burden of each AV and decreasing $\overline{\mathcal{P}}_1$. In contrast, fewer AVs lead to smaller interference, thus enhancing $\overline{\mathcal{P}}_2$. As a result, the overall performance $\overline{\mathcal{P}}_{\mathrm{overall}}$ can be enhanced under careful settings of $\lambda_a^0$ and $\check{d}$: i) given a $\check{d}$, an optimal value of $\lambda_a^0$, denoted by $\lambda_a^*$ (when $\overline{\mathcal{P}}_1=\overline{\mathcal{P}}_2$), can be used for the maximum $\overline{\mathcal{P}}_{\mathrm{overall}}$; ii) given $\lambda_a^0<\lambda_a^*$, the smaller $\check{d}$ can enhance $\overline{\mathcal{P}}_{\mathrm{overall}}$; and iii) given $\lambda_a^0>\lambda_a^*$, the larger $\check{d}$ can enhance $\overline{\mathcal{P}}_{\mathrm{overall}}$. 

Fig.~\ref{Fig-Ps}\subref{Fig-Ps-lambda-u}% and Fig.~\ref{Fig-Ps}
-\subref{Fig-Ps-theta} show the impact of $\lambda_u,\theta_i (\forall i=\{1, 2\}),\check{d}$ on $\overline{\mathcal{P}}_1, \overline{\mathcal{P}}_2, 
\overline{\mathcal{P}}_{\mathrm{overall}}$. It is observed that 
the increasing $\lambda_u$ reduces $\overline{\mathcal{P}}_1,\overline{\mathcal{P}}_{\mathrm{overall}}$ while it doesn't affect $\overline{\mathcal{P}}_2$; the increasing $\theta_i$ results in a decrease of $\overline{\mathcal{P}}_1,\overline{\mathcal{P}}_2,\overline{\mathcal{P}}_{\mathrm{overall}}$. This is because a larger $\lambda_u$ means more GUs, bringing more interference in first-hop transmissions, thereby reducing $\overline{\mathcal{P}}_1$. The second hop performance $\overline{\mathcal{P}}_2$ is related to the AV relays, then it won't be affected by $\lambda_u$. In summary, as $\lambda_u$ increases, $\overline{\mathcal{P}}_2$ decreases, and $\overline{\mathcal{P}}_2$ stay unchanged, then $\overline{\mathcal{P}}_{\mathrm{overall}}$ decreases. Meanwhile, a larger $\theta_i$ increases the difficulties of decoding signals in two-hop links, thus lowering $\overline{\mathcal{P}}_1,\overline{\mathcal{P}}_2$ and $\overline{\mathcal{P}}_{\mathrm{overall}}$. In Fig.~\ref{Fig-Ps}\subref{Fig-Ps-lambda-u}-\subref{Fig-Ps-theta}, $\overline{\mathcal{P}}_{\mathrm{overall}}$ increases with the growth of $\check{d}$. This is due to the relatively large value of $\lambda_a^0>\lambda_a^*$, as we observed in Fig.~\ref{Fig-Ps}\subref{Fig-Ps-lambda-a0}.
%The effect of $\check{d}$ on $\overline{\mathcal{P}}_1, \overline{\mathcal{P}}_2$ can refer to the observations in Fig.~\ref{Fig-Ps}\subref{Fig-Ps-lambda-a0}.

% From this figure, we have the following observations:
% \begin{itemize} 1) Given $\check{d}$ (i.e., the same line type), as $\theta_i$ increases, $\overline{\mathcal{P}}_1$ and $\overline{\mathcal{P}}_2$ gradually decreases, hence $\overline{\mathcal{P}}_{\mathrm{overall}}$ also decreases. This is because an increase in $\theta_i$ make it more difficult to decode signals from the GUs and AVs, respectively;
% 2) Given $\theta_i$, the lager $\check{d}$, the smaller $\overline{\mathcal{P}}_1$, the larger $\overline{\mathcal{P}}_2$, so as to $\overline{\mathcal{P}}_{\mathrm{overall}}$. The reasons are similar to those mentioned above. 
% \end{itemize}

%待定，午饭后修改
Overall, all simulation results closely align with the theoretical ones, validating the accuracy of our analytical model.
%\footnote{\color{black} To provide a better visualization, we present the simulation results of the case $\check{d}$ = 0 m, those of the cases $\check{d}$ = 100/200 m are similar to be obtained.}. 
Meanwhile, our observations offer essential insights for practitioners, such as network operators or engineers, in finding the optimal parameter settings to enhance system performance. For example, the observations in Fig.~\ref{Fig-Ps}\subref{Fig-Ps-lambda-a0} can assist in configuring AVs' deployments, including $\lambda_a^0$ and $\check{d}$, and then enhancing the overall two-hop transmission performance. 
%Based on the above observations, %we offer the following two \textbf{\textit{insights}}. \textit{First}, MHCPP is accurate and efficient in capturing the spatial separation attributes of AVs, which can assist in network modeling and practical deployments in GASS.
%\textit{Second}, the proposed theoretical model is helpful for practitioners, such as network operators or engineers, in finding the optimal parameter settings to enhance GASS performance.

\section{Conclusion}\label{Sec-Conclusion}
%This paper presents an overall connectivity analytical model to characterize the average performance of two-hop transmissions in GASS. Specifically, we adopt HPPP and MHCPP to model the location distributions of GUs and AVs. The GUs' associated distributions are explored under AVs, and the spherical locations of AVs are modeled to cater to practical network deployment. By applying stochastic geometry tools, we derive mathematical expressions of the spherical distance distributions for the two-hop transmission links, as well as the ASPs and overall connectivity of these links. Extensive Monte Carlo simulations have been conducted, and the results validate the accuracy of the proposed analytical model. %Our future work aims to explore more fine-grained performance metrics to enhance the understanding on overall network performance.

This paper presents an overall connectivity analytical model to characterize the average performance of two-hop transmissions in GASS. Specifically, we use HPPP and MHCPP to model the location distributions of GUs and AVs. We explore the distribution of GUs associated with AVs and model the spherical locations of AVs to reflect practical network deployment. By applying stochastic geometry tools, we derive mathematical expressions for the spherical distance distributions of the two-hop transmission links, as well as ASPs and the overall connectivity for these links. Extensive Monte Carlo simulations have been conducted, and the results validate the accuracy of the proposed analytical model. 

\section{Acknowledgment}
%{
%\small
The work %of Yalin Liu this paper 
was fully supported by the Hong Kong UGC/FDS project 
%UGC Faculty Development Scheme project 
under reference No. UGC/FDS16/E15/24. This work of Yulei Wang was completed when he was %affiliated 
at Hong Kong Metropolitan University. %(\textit{Corresponding author: Yalin Liu})
%}
\vspace{-0.3cm}

\section*{Appendix}
\label{app: angleproof1}
\textit{The proof of Proposition}~\ref{pro1}: 
\begin{align}
& \mathcal{P}_i
= \resizebox{0.9\hsize}{!}{$\mathbb{P}\left(\frac{P_i G_i h_i L_i(\texttt{t}_0)}{I_i +\sigma_i^2} > \theta_i \middle| \Phi_i\right) = 1- \mathbb{P}\left(h_i \le \frac{\theta_i \left(I_i +\sigma_i^2\right)}{P_i G_i L_i(\texttt{t}_0)} \middle| \Phi_i\right)$}
\notag\\
& \overset{(a)}{=}\resizebox{0.9\hsize}{!}{$ 1- \mathbb{E}_{I_i}\left[\frac{\gamma\left(m_i, \frac{m_i \theta_i \left(I_i + \sigma_i^2\right)}{\Omega_i P_i G_i L_i(\texttt{t}_0) } \right)}{\Gamma\left(m_i\right)} \right] 
\overset{(b)}{=} \mathbb{E}_{I_i}\left[\frac{ \Gamma\left(m_i, s_i \left(I_i +\sigma_i^2\right) \right) }{\Gamma\left(m_i\right)} \right] $} \notag\\
&\overset{(c)}{=} \mathbb{E}_{I_i} \left[ \exp \left(-s_i \left(I_i +\sigma_i ^2\right)\right) \sum_{k=0}^{m_i-1} \frac{\left(s_i \left(I_i + \sigma_i^2\right)\right)^{k}}{k!} \right] \notag\\
&\overset{(d)}{=} \sum_{k=0}^{m_i-1} {\frac{\left(-s_i\right)^k}{k!} \mathbb{E}_{I_i}\left[\frac{\mathrm{d}^k \exp\left(-s_i\left(I_i + \sigma_i^2\right)\right)}{\mathrm{d} s_i^k} \right]} \notag\\
&= \sum_{k=0}^{m_i -1} {\frac{\left(-s_i\right)^k}{k!} \left[\exp\left(-s_i\sigma_i^2\right) \mathcal{L}_{I_i}\left(s_i\right) \right]_{s_i}^{\left(k\right)}},
\label{Eq-PPP-STP2}
\end{align}%
where $s_i = m_i \theta_i/(\Omega_i P_i G_i L_i(\texttt{t}_0))$. %~\eqref{Eq-PPP-STP2}
$(a)$ arises from the cumulative distribution function of the gamma-distributed $h_i$. %~\eqref{Eq-PPP-STP2}
$(b)$ comes from the substitution of $\Gamma\left(s\right)=\gamma\left(s,x\right)+\Gamma\left(s,x\right)$, where $\gamma\left(s,x\right)$ and $\Gamma\left(s,x\right)$ are the lower and upper incomplete gamma function, respectively. %~\eqref{Eq-PPP-STP2}
$(c)$ is resulted from the incomplete gamma function~\cite[8.352.2]{gradshteyn2014table}. %$\Gamma(1+n, x) = n! e^{-x} \sum_{m=0}^{n} \frac{x^{m}}{m!}$ in \cite{gradshteyn2014table}. 
%~\eqref{Eq-PPP-STP2}
$(d)$ arises from the substitution of $\exp\left(-sx\right)x^k=\left(-1\right)^k\frac{\mathrm{d}^k \exp\left(-sx\right)}{\mathrm{d} s^k}$~\cite{galkin2019stochastic}.  
% $(e)$ comes from the Leibniz integral rule. $(e)$ is from fact that the  $k$-th derivative of $\mathcal{L}_{(I+\sigma^2)}\left(s\right)$ can be expressed as a sum of products of higher derivatives of $\mathcal{L}_I\left(s\right)$ and $\exp(-s\sigma^2)$  following the general Leibniz rule \cite{galkin2019stochastic},\cite{galkin2021reqiba}.
In~\eqref{Eq-PPP-STP2}, $\mathcal{L}_{I_i}\left(s_i\right)$ is the Laplace transform of the interference $I_i$, which can be further calculated by
\begin{align}
& \mathcal{L}_{I_i}\left(s_i \right) =\mathbb{E}_{I_i}\left[\exp\left(-s_i I_i\right)\right] \notag \\
%&=\mathbb{E}_{r_i(\texttt{t}),h_i}\left[\exp\left(-s_i \sum_{\texttt{t} \in \Phi_i^\prime}{P_i G_i h_i L_i(\texttt{t})} \right)\right] \notag \\
&\overset{(a)}{=} \mathbb{E}_{r_i(\texttt{t})}\left[\mathbb{E}_{h_i}\left[ \exp \left(-\sum_{\texttt{t} \in \Phi_i^\prime}{s_i P_i G_i h_i L_i(\texttt{t}) }\right)\right]\right] \notag  \\ 
&\overset{(b)}{=} \mathbb{E}_{r_i(\texttt{t})}\left[\prod_{\texttt{t} \in \Phi_i^\prime}{\mathbb{E}_{h_i}\left(\exp \left(-s_i P_i G_i L_i(\texttt{t}) h_i \right) \right) }\right] \notag \\
&\overset{(c)}{=} \mathbb{E}_{r_i(\texttt{t})}\left[\prod_{\texttt{t} \in \Phi_i^\prime} \left(1 + s_i Q_i L_i(\texttt{t}) \right)^{-m_i}  \right],
\label{Eq-LT-1}
\end{align}%
where $ Q_i = P_i G_i \Omega_i /(m_i)$. %~\eqref{Eq-LT-1}
$(a)$ comes from the fact that $r_i(\texttt{t})$ and $h_i$ are mutually independent. %~\eqref{Eq-LT-1}
$(b)$ follows from the property of exponential distribution, i.e., $\exp\left(\sum_{j} h_i\right)=\prod_{j} \exp\left(h_i\right)$. %~\eqref{Eq-LT-1}
$(c)$ comes from the moment-generating function of $h_i$~\cite{simon2004digital}. %, i.e., $\mathbb{E}[\exp(sH)=(1-\frac{s\Omega}{m})^{-m}]$ \cite{simon2004digital}. 
According to the probability generating functional of an HPPP, % $\Phi$ with density $\lambda$, 
i.e., $\mathbb{E}\left[\prod_{x\in\Phi} f\left(x\right)\right]=\exp\left(-\lambda\int_{\mathbb{R}^d}\left(1-f\left(x\right)\right)\mathrm{d} x\right)$, 
Eq.~\eqref{Eq-LT-1} can be calculated by %we express the last equation of~\eqref{Eq-LT-1} as
\begin{equation}
\begin{aligned}
& \mathcal{L}_{I_i}\left(s_i\right) 
\overset{(a)}{=} \resizebox{0.8\hsize}{!}{$ \exp \left(-\lambda_i \int_{\mathcal{A}_i }{\left(1-\left(1 + s_i Q_i L_i(\texttt{t}) \right)^{-m_i} \right) \mathrm{d} \texttt{t} }\right) $} \\
& \overset{(b)}{=} \resizebox{0.95\hsize}{!}{$ \exp \left(-\lambda_i  \int_{0}^{2\pi }{\int_{0}^{\varphi_{\texttt{t}} }{\int_{\check{R}}^{\hat{R}}{\left(1- \left(1 + s_i Q_i L_i(\texttt{t}) \right)^{-m_i} \right) r^2 \mathrm{sin}\varphi \mathrm{d}r} \mathrm{d}\varphi } \mathrm{d}\vartheta }\right) $}  \\
& \overset{(c)}{=} \resizebox{0.95\hsize}{!}{$ \exp \left(-2 \pi \lambda_i \Tilde{R}_{\mathrm{e}}^2 \int_{0}^{\varphi_{\texttt{t}}}{\left(1- \left( 1 + s_i Q_i L_i(\texttt{t}) \right)^{-m_i} \right) \mathrm{sin} \varphi \mathrm{d}\varphi }\right).  $}
\label{Eq-LT-2}
\end{aligned}
\end{equation}
In~\eqref{Eq-LT-2}, $(a)$ comes from the integral of the interfering GU (or AV) $\texttt{t}$ in $\mathcal{A}_i$. %$\mathcal{A}_i$ is the area in which the interfering GUs (or AVs) locate; $\lambda_i$ is the density of interfering GUs (or AVs), $\lambda_1 = \lambda_u^t$, and $\lambda_2 = \lambda_a$; 
%~\eqref{Eq-LT-2}
$(b)$ is resulted from the conversion from Cartesian coordinates $\texttt{t}: (x, y, z)$ to polar coordinates $\texttt{t}: (r, \varphi, \vartheta)$, i.e., $\mathrm{d} \texttt{t} = \mathrm{d}x \mathrm{d}y \mathrm{d}z = r^2 \mathrm{sin}\varphi \mathrm{d}r \mathrm{d}\varphi \mathrm{d}\vartheta$, $r \in [\check{R}, \hat{R}]$ is the radial distance, $\vartheta \in [0, 2\pi]$ is the azimuthal angle, $\varphi \in [0, \varphi_{\texttt{t}}]$ is the polar angle. In %~\eqref{Eq-LT-2}
$(c)$, $\check{R} = \hat{R} = \Tilde{R}_{\mathrm{e}}$. Because the GUs/AVs located on the spherical surface with radius $\Tilde{R}_{\mathrm{e}}$, and $$\Tilde{R}_{\mathrm{e}} = \left\{\begin{matrix} R_{\mathrm{e}}, & i = 1 \\ R_{\mathrm{e}} + H_a, & i = 2 \\  \end{matrix}.\right.$$%
Substituting~\eqref{Eq-LT-2} into~\eqref{Eq-PPP-STP2} and following the similar derivation process in~\cite{liu2024space}, we %derive the final expression of $\mathcal{P}_{i}\left(\theta_i\right)$ in
have~\cref{pro1}. %Due to page limitations, we put the detailed derivation process into supplementary material~\cite{VTC25Supplementary}.
\hfill$\blacksquare$	

% References
\bibliographystyle{IEEEtran}
\bibliography{reference.bib}

\end{document}